\documentclass[12pt,amsymb,fullpage]{amsart}
\usepackage{amssymb,amscd,pstricks}

\newtheorem{theorem}{Theorem}[section]
\newtheorem{defn}[theorem]{Definition}

\newtheorem{lemma}[theorem]{Lemma}
\newtheorem{fact}[theorem]{Fact}

\newtheorem{eple}[theorem]{Example}
\newtheorem{rmk}[theorem]{Remarks}
\newtheorem{dsc}[theorem]{Discussion}
\newtheorem{nota}[theorem]{Notation}

\newsavebox{\indbin}
\savebox{\indbin}{\begin{picture}(0,0)
\newlength{\gnu}
\settowidth{\gnu}{$\smile$} \setlength{\unitlength}{.5\gnu}
\put(-1,-.65){$\smile$} \put(-.25,.1){$|$}
\end{picture}}

\newcommand{\be}{\begin{enumerate}}
\newcommand{\bd}{\begin{defn}}
\newcommand{\bt}{\begin{theorem}}
\newcommand{\bl}{\begin{lemma}}
\newcommand{\ee}{\end{enumerate}}
\newcommand{\ed}{\end{defn}}
\newcommand{\et}{\end{theorem}}
\newcommand{\el}{\end{lemma}}

\begin{document}
\title{Zariski Structures and Algebraic Geometry}
\author{Tristram de Piro}
\address{Mathematics Department, James Clerk Maxwell Building, Kings Buildings, Edinburgh, }
\email{depiro@maths.ed.ac.uk}
\thanks{The author was supported by the Seggie Brown Research Fellowship}
\begin{abstract}
The purpose of this paper is to provide a new account of
multiplicity for finite morphisms between smooth projective
varieties. Traditionally, this has been defined using commutative
algebra in terms of the length of integral ring extensions. In
model theory, a different approach to multiplicity was developed
by Zilber using the techniques of non-standard analysis. Here, we
first reformulate Zilber's method in the language of algebraic
geometry using specialisations and secondly show that, in
classical projective situations, the two notions essentially
coincide. As a consequence, we can recover intersection theory in
all characteristics from the non-standard method and sketch the
further development of the theory in connection with etale
cohomology  and deformation theory. The usefulness of this
approach can be seen from the increasing interplay between Zariski
structures and objects of non-commutative geometry, see
\cite{Zilb1}.
\end{abstract}

\maketitle

We will work mainly in the language of Weil's Foundations, namely
using varieties instead of schemes. $K$ will denote a big
algebraically closed field. $L\subset K$ will denote a small
algebraically closed field. By an affine variety $V$, we mean a
closed subset of $K^{n}$ in the Zariski topology. If $V$ is
irreducible, we denote the ring of regular functions on $V$ by
$K[V]$ and the function field by $K(V)$. If $k\subset K$ is
perfect, we say that $V$ is defined over $k$ if $I(V)$, the
radical ideal of functions vanishing on $V$ is generated by
polynomials with coefficients in $k$. Any irreducible affine
variety $V$ has a minimal field of definition $k_{V}$ with the
property that any automorphism fixes $V$ setwise iff it fixes
$k_{V}$ pointwise. This is a classical result due to Weil, but is
in fact a special case of a more general construction due to model
theorists of canonical bases, see \cite{Kim}. By a variety, we
will mean a set $V$, a covering of subsets $V_{1},\ldots V_{m}$
and for each $i$ a bijection $f_{i}:V_{i}\rightarrow U_{i}$ with
$U_{i}$ an affine variety and such that for each $1\leq i,j\leq
m$, $U_{ij}=f_{i}(V_{i}\cap V_{j})$ is an open subset of $U_{i}$
and $f_{ij}=f_{j}f_{i}^{-1}$ is an isomorphism between the affine
varieties $U_{ij}$ and $U_{ji}$.  A variety $V$ then inherits a
natural Zariski topology by declaring $U\subset V$ open if for
each $i$, $f_{i}(U\cap V_{i})$ is open in $U_{i}$. For $k\subset
K$, we will say that $V$ is defined over $k$, if the data
$(U_{i},U_{ij},f_{ij})$ is defined over $k$ in the sense of affine
varieties. We let $P^{n}(K)$ denote $n$-dimensional projective
space over $K$, that is $K^{n+1}/{\verb1~1}$, where ${\verb1~1}$
is the equivalence relation on $K^{n+1}\setminus \{\bar 0\}$ given
by $(x_{0},\ldots, x_{n}){\verb1~1}(y_{0},\ldots y_{n})$ iff
$\lambda(x_{0},\ldots, x_{n})=(y_{0},\ldots,y_{n})$ for some
$\lambda\in K$. Writing elements of $P^{n}(K)$ in homogenous
coordinates, $(x_{0}:x_{1}:\ldots:x_{n})$, we have natural
bijections $f_{i}$ between $K^{n}$ and $P^{n}(K)_{i}=\{\bar
x:x_{i}\neq 0\}$. This gives $P^{n}(K)$ the structure of a variety
defined over the prime subfield and an induced Zariski toplogy. By
a projective variety $V$, we mean a closed subset of $P^{n}(K)$,
using the coordinate charts $f_{i}$, $V$ automatically is a
variety in the sense defined above. Equivalently, a projective
variety $V$ is defined by a set of homogenous polynomials in
$K[x_{0},\ldots x_{n}]$ and is defined over $k$ if the ideal
$I(V)$ is generated by homogenous polynomials with coefficients in
$k$. If a variety is $V$ defined over $k$ and $k\subset L\subset
K$ with $L$ algebraically closed then we will use the notation
 $V(L)$ to denote $V$ considered as a variety over $L$. In this case, we will
 require that a subvariety of $V(L)$ is  defined over $L$.\\
   We will use the notation $X\times_{Y}Z$ to denote the fibre product of two varieties $X$ and $Y$ over $Z$. Given a variety $V$ defined over $k$ and a tuple of elements $\bar a\in V^{n}$, we will use $k(\bar a)$ to denote the field of definition of $\bar a$. In the case when $X=Spec(L)$, corresponding to an $L$ rational point $j:Spec(L)\rightarrow Z$;\\

\newcommand{\End}{\operatorname{End}}
\begin{eqnarray*}
\begin{CD}\\
L\times_{Z}Y@>i>>X\\
   @VVprV   @VVjV\\
Y@>f>>Z\\
\end{CD}
\end{eqnarray*}\\

we will often use the notation $L\times_{Z}Y$ to denote the
geometric fibre $f^{-1}(y)$ of a point $y\in Z$, considered as a
variety over $L$. Similar notation will be used in the case of
sheaves. Given varieties $Y$,$Z$ and a morphism $g:Y\rightarrow
Z$,  we define the pullback of a coherent sheaf $F$ on $Z$ to be
the sheafification of \\

$g^{*}F=O_{Y}\otimes_{g^{-1}O_{Z}}g^{-1}F$\\

where $g^{-1}F(U)=lim_{\rightarrow, g(U)\subset V}F(V)$. Again, in
the case when $j:Spec(L)\rightarrow Z$ is an $L$ rational point
and $F$ is a coherent sheaf on $Z$, $j^{-1}F=F_{z}$,the localised
sheaf at $z$, and $L\otimes_{O_{z,Z}}F_{z}$ is a vector space over
$L$ which, by slight abuse of notation, corresponds to the fibre
of the sheaf $F$ at $z$. Given a morphism $f:X\rightarrow Y$, we
let $\Omega_{X/Y}$ denote the sheaf of relative differentials on
$X$. We will use the geometric construction of $\Omega_{X/Y}$ as
$\Delta^{*}J/J^{2}$ where $\Delta:X\rightarrow X\times_{Y}X$ is
the diagonal embedding and $J/J^{2}$ is the normal bundle of
$\Delta(X)$ in $X\times_{Y}X$. In the case when $Y=Spec(L)$ for
$k\subset L\subset K$ and $k$ the field of definition of $X$, we
use the notation $\Omega_{X/L}$ to denote the sheaf of meromorphic
differentials on $X$ and $\Omega_{X/L}^{*}$ the sheaf of
meromorphic vector fields. There is a canonical isomorphism;\\

$d:m_{z}/m_{z}^{2}\rightarrow (\Omega_{X/L})_{z}\otimes L$\\

$d(f+m_{z}^{2})=df$\\

relating the sheaf of differentials to the cotangent space at a point. Using this isomorphism and Nakayama's Lemma, one has that for an algebraic variety $X$ of dimension $n$ over $k\subset L$, $\Omega_{X/L}$ is a locally free module of rank $n$ on the nonsingular locus $U$ of $X$, see \cite{Mum} for details.\\

\begin{section}{Zariski Structures}

\begin{defn}

Let $({\mathcal M},\tau)$ be a topological space and let $\{C\}$ denote the collection of closed sets on $({\mathcal M})$. We call $(\mathcal M,\tau)$ a Zariski structure if the following axioms hold;\\

\end{defn}

(L) Language: Basic relations are closed;\\

The diagonals $\Delta_{i}\subset {\mathcal M}^{i}\times {\mathcal M}^{i}$ are closed.\\

Any singleton in ${\mathcal M}$ is closed.\\

Cartesian products of closed sets are closed\\

(P) Properness: The projection maps $pr:{\mathcal M}^{n+1}\rightarrow {\mathcal M}^{n}$ are proper and continuous, that is the images and inverse images of closed sets under $pr$ are closed\\

(DCC) Descending Chain Condition: The topology given by the closed sets on ${\mathcal M}^{n}$ is Noetherian for all $n\geq 1$. The condition $(DCC)$ implies that every closed set $C$ can be written uniquely (up to permutation) as a union of irreducible closed sets;\\

$C=C_{1}\cup\ldots\cup C_{n}$\\

(DIM) Dimension: The following notion of dimension for closed sets $C\subset {\mathcal M}^{n}$ is well defined;\\

For irreducible $C$, $dim(C)$ is the maximum $m$ for which there exists a chain of irreducible closed sets $C_{0}\subset C_{1}\subset\ldots C_{m}=C$.\\

For arbitrary closed $C$, $dim(C)=max_{1\leq i\leq m}\{dim(C_{i})\}$ for $C_{i}$ the irreducible components of $C$\\

(PS) Pre-Smoothness: For all closed irreducible sets $C_{1},C_{2}\subset {\mathcal M}^{n}$, with $C_{1}\cap C_{2}\neq\emptyset$, \\

$dim_{comp}(C_{1}\cap C_{2})\geq dim(C_{1})+dim(C_{2})-dim({\mathcal M}^{n})$\\

(DF) Definability of fibres: If $C\subset{\mathcal M}^{n+m}$ is closed, then\\

$F(C,k)=\{\bar a\in{\mathcal M}^{n}:dim(C(\bar a))>k\}$\\

is closed.\\

(GF) Generic fibres: If $C\subset {\mathcal M}^{n+m}$ is closed and irreducible, then\\

$dim(C)=dim(pr(C))+min_{\bar a\in pr(C)}dim C(\bar a)$\\

\begin{rmk}

The definition of dimension easily implies the following properties;\\

(DU) Dimension of unions: For $C_{1},C_{2}$ closed, then\\

$dim(C_{1}\cup C_{2})=max\{dim(C_{1}), dim(C_{2})\}$\\

(DP) The dimension of a point is $0$.\\

(DI) Dimension of irreducible sets: If $C_{1}\subsetneq C_{2}$ and $C_{2}$ is irreducible, then $dim(C_{1})<dim(C_{2})$.\\

\end{rmk}

We now show the following;\\

\begin{theorem}

Let $V$ be a smooth projective variety of dimension $m$ defined
over $k$ and $k\subset L$ with $L$ algebraically closed, then
$V(L)$ considered as a toplogical space with closed sets given by
the algebraic subvarieties defined over $L$ is an irreducible
Zariski structure of dimension $m$.

\end{theorem}

\begin{proof}

We will verify the axioms;\\

(L) We need only verify that the diagonals $\Delta_{i}\subset V^{i}\times V^{i}$ are closed. \\

(P) An algebraic variety $V$ is complete if for all varieties $Y$, the projection morphism\\

$pr: V\times Y\rightarrow Y$\\

is closed. Taking $Y$ to be $V^{n}$ in the above definition, complete varieties have the property that the projection maps\\

$pr: V^{n+1}\rightarrow V^{n}$\\

are closed. If $W\subset V$ is a closed subvariety of a complete variety $V$, then, as is easily checked, $W$ is also complete. By assumption $V$ is a closed subvariety of $P^{N}(L)$ for some $N$. Now it is a classical fact that $P^{N}(L)$ is complete, see for example \cite{Humph}.\\

(DCC) Let $\{W_{i}:i<\omega\}$ be an infinite descending chain of
closed
subvarieties of $V^{n}$. Let $\{U_{1},\ldots U_{n}\}$ be an affine open cover of $V^{n}$. Then $\{U_{j}\cap W_{i}:i<\omega\}$ defines a descending chain of closed subvarieties of each $U_{j}$. By the Nullstellensatz, each such chain stabilises inside $U_{j}$. Then clearly the chain stabilises inside $V^{n}$. \\

(DIM) For $W$ an irreducible subvariety of $V^{n}$, we let
$dim_{geom}(W)=t.deg(L(W)/L)$. Then $dim_{geom}$ corresponds to
$dim$ as defined above. To see this, suppose that $dim(W)\geq
n+1$, and $W$ is irreducible, then by definition one can find an
irreducible closed subvariety $W'\subset W$ with $dim(W')\geq n$
and so inductively $dim(W')\geq n$. Now take any affine open
subset of $V^{n}$ intersecting $W'$, so we may assume that $W$ and
$W'$ are affine as the function fields are unchanged. Let $L[W]$
denote the coordinate domain of $W$, $p$ the proper prime ideal
corresponding to $W'$ and $dim_{Krull}$ the Krull dimension of an
integral domain. By Krull's theorem,
$height(p)+dim_{Krull}(L[W]/p)=dim_{Krull}(L[W])$,
$dim_{Krull}(L[W])=t.deg(L(W))$ and
$dim_{Krull}(L[W]/p)=t.deg(L(W'))$, hence
$t.deg(L(W')))<t.deg(L(W))$. It follows that
$dim_{geom}(W')<dim_{geom}(W)$ and so $dim_{geom}(W)\geq n+1$.
Conversely, if $dim_{geom}(W)\geq n+1$, then again assuming $W$ is
irreducible and affine, if we take $f\in L[W]$ to be a non-unit, then each irreducible component of $V(f)\subset W$ has codimension $1$ in $X$, see \cite{Humph}. Therefore, $dim_{geom}(V(f))\geq n$ and inductively $dim(V(f))\geq n$. As each component of $V(f)$ is a proper closed subset of $X$, $dim(W)\geq n+1$. Now clearly we have that $dim_{geom}$ corresponds to $dim$ and so in particular we know that $dim(V^{n})=mn$ and the notion of $dim$ on $V^{n}$ is well defined.\\

(PS) A simple calculation shows that for $(x_{1}\ldots x_{n})\in V^{n}$, $m_{\bar x}\cong \Sigma_{i=1}^{n}O_{x_{1}\ldots \hat {x_{i}}\ldots x_{n}}\otimes m_{x_{i}}$. Hence, \\

$Tan_{\bar x}(V^{n})=(m_{\bar x}/m_{\bar
x}^{2})^{*}\cong\Sigma_{i=1}^{n}(m_{x_{i}}/m_{x_{i}}^{2})^{*}=\Sigma_{i=1}^{n}Tan_{x_{i}}(V)$.

Therefore, $V^{n}$ is smooth.\\

Now we use the following lemma;\\

\begin{lemma}

If $X$ is a non-singular algebraic variety of dimension $n$, and $Y,Z$ are irreducible closed subsets. Then if $W$ is a component of $Y\cap Z$, we have,\\

$dim(W)\geq dim(Y)+dim(Z)-n$\\

or equivalently\\

$codim(W)\leq codim(Y)+codim(Z)$\\

\end{lemma}

\begin{proof}
We have that $Y\cap Z\cong Y\times Z\cap\Delta(X)$ inside
$X\times_{L}X$. Let $g_{1},\ldots ,g_{n}$ be uniformizers on an
open subset $U$ inside $X$. Then we saw above that $\Omega_{X/L}$
is just the pullback of the conormal sheaf $J/J^{2}$ for the
inclusion of $\Delta(X)$ inside $X\times_{L} X$. As $\Omega_{X/L}$
is locally free, so is $J/J^{2}$, and in particular generated
freely on $\Delta(U)$ by the functions $g_{1}\otimes 1-1\otimes
g_{1},\ldots, g_{n}\otimes 1-1\otimes g_{n}$. At a point
$x\in\Delta(U)$, we have that $g_{1}\otimes 1-1\otimes
g_{1},\ldots, g_{n}\otimes 1-1\otimes g_{n}$ generate
$J_{x}/J_{x}^{2}$ and therefore form a basis for the vector space
$J_{x}/m_{x}J_{x}$ as clearly any function belonging to $J_{x}$
lies in $m_{x}$ the ideal of functions in $O_{X\times X,x}$
vanishing at $x$. Then, as $J_{x}/m_{x}J_{x}$ is just the base
change $J\otimes k(x)$ of the ideal sheaf $J$ at the point $x$, it
follows by Nakayama's lemma that these functions generate $J$ on an open neighborhood $U$
 containing $x$ (not freely!). It follows that $Y\times Z\cap\Delta(X)$ is cut out by exactly
 $n$ equations inside $Y\times Z$, so by standard dimension theory we have the result.\\

\end{proof}

It follows immediately that $V^{n}$ satisfies $(PS)$.\\

In order to check the final $2$ axioms we introduce the following definitions;\\

\begin{defn}

If $\bar a,\bar b\in V^{n}$ are tuples of elements, we define
$locus(\bar a/\bar b)$ to be the intersection of all closed
subvarieties defined over $k(\bar b)$ containing $\bar a$ and
$locus_{irr}(\bar a/\bar b)$ to be the intersection of all closed
subvarieties defined over $k(\bar b)^{alg}$. We define $dim(\bar
a/k)$ to be $t.deg(k(\bar a)/k)$ and $dim(\bar a/k\bar b)$ to be
$t.deg(k(\bar a)/k(\bar b))$; if the underlying field $k$ is clear
from context, we will abbreviate this to $dim(\bar a/\bar b)$

\end{defn}

By the condition $DCC$, it is clear that $locus$ and $locus_{irr}$ are well defined. $locus_{irr}$ is an irreducible subvariety of $V^{n}$ containing $\bar a$, as if $V$ is an irreducible components of $locus_{irr}$ containing $\bar a$ and $k_{V}$ is the minimal field of definition, then $k_{V}$ has only finitely many congugates under an automorphism fixing $k(\bar b)^{alg}$, hence $k_{V}\subset k(\bar b)^{alg}$.\\

\begin{defn}

If $W\subset V^{n}$ is an irreducible closed subvariety, $\bar
a\in W$, and $\bar b$ is a tuple of elements such that $k(\bar b)$
contains a field of definition for $W$, then we say that $\bar a$
is generic in $W$ over $\bar b$ if $locus(\bar a/\bar b)=W$.

\end{defn}

\begin{lemma}

Let $W\subset V^{n}$ be an irreducible closed subvariety defined
over $k$ and $\bar a$ generic in $W$ over $k$. Then
$dim(W)=t.deg(k(\bar a)/k)=dim(\bar a/k)$.

\end{lemma}

\begin{proof}

By the above, $dim(W)=dim_{geom}(W)=t.deg(k(W)/k)$. By choosing an
open affine subvariety of $W$ containing $\bar a$ and defined over
$k$, we can assume that $W$ is affine. Now define a map
$ev:k[W]\rightarrow k(\bar a)$ by setting $ev(f)=f(\bar a)$. $ev$
is injective as if $f(\bar a)=0$, then as $f$ has coefficients in
$k$ and $\bar a$ is generic in $W$ over $k$, $f|W=0$. Clearly $ev$
extends to a map on $k(W)$ which is an isomorphism.

\end{proof}

(DF)  Let $W\subset V^{n+m}$ be a closed subvariety and $pr$ the projection onto $n$ factors. We can cover $(P^{N}(L))^{(n+m)}$ with finitely many affines of the form $A^{N(n+m)}$, hence we may assume that $W$ is a closed subvariety of $A^{N(n+m)}$ and show that $\Gamma(\bar y)=\{\bar a: dim(W(\bar a))\geq k+1\}$ is closed in $pr(W)$. By additivity of $t.deg$ and the lemma, this occurs iff we can find algebraically independent elements $b_{1}\ldots b_{k}b_{k+1}\subset \bar b\subset L$ such that $W(\bar b\bar a)$ holds iff  \\

$\exists_{\sigma(k+2)}\ldots\exists_{\sigma(Nm)} W(x_{1},\ldots,
x_{Nm},\bar a)
$\\

has maximal dimension for some permutation $\sigma\in S_{Nm-(k+1)}$. We may write each projection $W_{\sigma}$ in the form \\

$\bigcap_{i} F_{i}(x_{1}\ldots x_{k+1},\bar y)=0 \cap \bigcap_{j}Q_{j}(x_{1}\ldots x_{k+1},\bar y)\neq 0$\\

where $F_{i}$ and $Q_{j}$ are polynomials in the variables $\bar x\bar y$. Let $\theta_{\sigma}(\bar y)$ define the closed set given by the vanishing of all coefficients in the $F_{i}$. Then an easy calculation shows that $\Gamma_{\sigma}(\bar y)=\{\bar y\in pr(W):\theta_{\sigma}(\bar y)\}$, which is closed.\\

(GF) We first show the following;\\

\begin{lemma}

Let $W\subset V^{n+m}$ be closed and irreducible, defined over
$k$. Then $\bar a\bar b$ is generic in $W$ over $k$ iff $\bar a$
is generic in $pr(W)$ over $k$ and $\bar b$ is generic in $W(\bar
a)$ over $k(\bar a)$.

\end{lemma}

\begin{proof}
One direction is straightforward, if $\bar a$ is not generic in $pr(W)$,
then $\bar a\in E\subsetneq pr(W)$ and $\bar a\bar b\in pr^{-1}(E)\subsetneq W$.
 If $\bar b$ is not generic in $W(\bar a)$, then we can find $X\subsetneq W(\bar a)$
containing $\bar b$ defined over $k(\bar a)$. As we are working in a product of $P^{m}(L)$, we can
define $X$ by a series of $n$-homogeneous equations with coefficients in $k(\bar a)$. Applying
 Frobenius to these equations, we can in fact assume that the coefficients lie in $k<\bar a>$. Now a
 straightforward exercise in clearing denominators and writing affine equations in
 homogeneous form shows that we can write $X$ as the fibre $Y(\bar a)$ for some
 closed subvariety $Y$ of $V^{n+m}$. Intersecting with $W$ if necessary gives a proper
closed $Y\subsetneq W$ with $\bar a\bar b\in W$ and defined over $k$. \\

For the other direction, suppose that $\bar a\bar b$ is not generic in $W$ over $k$, then there exists $X$ defined over $k$ such that $\bar a\bar b\in X\subsetneq W$. Then $\bar a\in pr(X)$ which is also closed and defined over $k$. Hence, $pr(X)=pr(W)$. As $\bar b\in X(\bar a)$, we have that $dim(X(\bar a))=dim(W(\bar a))=m$. By $(DF)$,\\

$X_{m}=\{\bar a\in pr(X): dim(X(\bar a))=dim(W(\bar a))=m\}$\\

is constructible and, by automorphism, can be seen to be defined
over $k$. Hence, as $\bar a$ was assumed to be generic, $X_{m}$ is
open inside $pr(X)$. Now, using Lemma 1.7 and the hypotheses on
$\bar a$,$\bar b$, $dim(X)\geq dim(\bar a\bar b/k)=dim(\bar b/\bar
a k)+dim(\bar a/k)=m+dim(pr(W))$. However, choosing $\bar a'\bar
b'$ generic in $W$ over $k$, we have that $dim(W)=dim(\bar a'\bar
b'/k)=m+dim(pr(W))$ by the properties of $X_{k}$. Hence,
$dim(X)\geq dim(W)$ contradicting the fact that $X\subsetneq W$
and $W$ was assumed to be irreducible.

\end{proof}

Using the lemma, we can give an easy proof of $(GF)$;\\

Let $W\subset V^{n+m}$ be closed, irreducible and defined over $k$. Choose $\bar a\bar b$ generic in $W$ over $k$. Then \\

$dim(W)=dim(\bar a\bar b/k)=dim(\bar b/\bar a k)+ dim(\bar a/k)=dim(pr(W))+ min_{\bar a\in pr(W)}W(\bar a)$. \\

The last equality follows from the previous lemma and $(DF)$.\\

We have therefore checked all the axioms.

\end{proof}

\begin{defn}

Given a closed subvariety $W$ of $V^{m}(L)$ and a closed $F\subset
W\times V^{m}$, all defined over $k$, we say that $F$ is a cover
of $W$ if $pr(F)=W$ and that $\bar a\in W$ is regular for the
cover if $dim F(\bar a)=dim F(\bar a')$ for $\bar a'$ generic in
$W$ over $k$.

\end{defn}

\end{section}

\begin{section}{Specialisations}

In order to apply the technique of specialisations, we fix an algebraically closed field $L$ and construct a universal extension $K_{\omega}$ as follows.\\

Set $L=K_{0}$. Construct $K_{i+1}$ inductively as follows;\\

Let $K_{i}((t_{i+1}))$ be the field of formal Laurent series in the variable $t_{i+1}$ over the algebraically closed field $K_{i}$. Define $K_{i+1}=K_{i}((t_{i+1}))^{alg}$.\\

Given the tower of algebraically closed fields $L\subset K_{1}\subset K_{2}\subset\ldots\subset K_{i}\subset\ldots$, we set $K_{\omega}=\bigcup_{i<\omega}K_{i}$\\

For all $i<\omega$, $K_{i+1}$ is equipped with a canonical valuation $\overline{v_{i+1}}:K_{i+1}\rightarrow {\mathcal Z}$ defined as follows;\\

For $f\in K_{i}((t_{i+1}))$, we set $v_{i+1}(f)=ord_{i+1}(f)$, where $ord_{i+1}(f)$ is the minimum $n$ appearing in the Laurent expansion of $f$. As is shown in \cite{Neuk}, $(K_{i}((t_{i+1})), v_{i+1})$ is the completion of $(K_{i}(t_{i+1}), v)$, for the canonical valuation $v$ on the function field $K_{i}(t_{i+1})$. It follows that $K_{i}((t_{i+1}))$ is a Henselian field with respect to $v_{i+1}$. By Hensel's lemma, $K_{i}((t_{i+1}))^{alg}$ is a union $\bigcup_{i<\omega}K_{i}((t_{i+1}^{1/n}))$ of ramified extensions of $K_{i}((t_{i+1}))$. The valuation $v_{i+1}$ extends uniquely to the spectral valuation on $K_{i}((t_{i+1}^{1/n}))$ by the formula;\\

$\overline{v_{i+1}}(\alpha)=(1/n) v_{i+1}(N_{K_{i}((t^{1/n}))/K_{i}((t))}(\alpha))$\\

From the previous section, we have that $P^{n}(L)$ and $P^{n}(K_{\omega})$ with closed sets given by subvarieties defined over $L$ and $K_{\omega}$ respectively are Zariski structures. We define a specialisation map $\pi_{\omega}: P^{n}(K_{\omega})\rightarrow P^{n}(L)$ as follows. First, the maps;\\

$\pi_{i+1,m}:P^{n}(K_{i}((t_{i+1}^{1/m})))\rightarrow P^{n}(K_{i})$\\

are defined by\\

$(f_{0}:\ldots:f_{n})\mapsto (\overline {t^{s}f_{0}}:\ldots:\overline { t^{s}f_{n}})$\\

where $s\in{ \mathcal Z}$ is chosen such that $\{t^{s}f_{0},\ldots t^{s}f_{n}\}\subset O_{v_{i+1}}=K_{i}[[t_{i+1}^{1/m}]]$ and $v_{i+1}(t_{i+1}^{s}f_{j})=0$ for some $j$ with $0\leq j\leq n$. Using the fact that the residue mapping is a homomorphism on $K_{i}((t_{i+1}^{1/m}))$, this map is clearly well defined. Moreover, the maps $\pi_{i+1,m}$ are compatible for $m\in {\mathcal N}$, in the sense that, given $m_{1}$ and $m_{2}$, for $\pi_{i+1,m_{1}m_{2}}:P^{n}(K_{i}((t_{i+1}^{1/m_{1}m_{2}})))\rightarrow P^{n}(K_{i})$, we have that $\pi_{i+1,m_{1}m_{2}}|P^{n}(K_{i}(t_{i+1}^{1/m_{k}}))=\pi_{i+1,m_{k}}$ for $k\in\{1,2\}$. Hence, the maps $\pi_{i+1,m}$ naturally define a map \\

$\pi_{i+1}:P^{n}(K_{i+1})\rightarrow P^{n}(K_{i})$\\

$(f_{0}:\ldots:f_{n})\mapsto \pi_{i+1,m}(f_{0}:\ldots:f_{n})$\\

where $\{f_{0},\ldots,f_{n}\}\subset K_{i}[[t^{1/m}]]$\\

Now, for $M<N$, let  $\Pi_{N,M}=\pi_{M+1}\circ\ldots\circ\pi_{N}:P^{n}(K_{N})\rightarrow P^{n}(K_{M})$ and let $\Pi_{M,M}=Id$. Then we have that for $M_{1}\leq M_{2}\leq M_{3}$, $\Pi_{M_{1},M_{2}}\circ \Pi_{M_{2},M_{3}}=\Pi_{M_{1},M_{3}}$, hence the maps $\{\Pi_{N,M}\}$ form an inverse system and we set $\Pi_{\omega}=\bigcup_{M,N}\Pi_{M,N}:P^{n}(K_{\omega})\rightarrow P^{n}(L)$.\\

We now show the following lemmas for the pair $(P^{n}(K_{\omega}),\Pi_{\omega})$. \\

\begin{lemma}

Let $V\subset P^{n}(L)$ be a smooth projective variety defined
over $L$. Then $\Pi_{\omega}:V(K_{\omega})\rightarrow V(L)$
defines a homomorphism of Zariski structures, in the sense that
for all closed $W\subset V^{m}$ defined over $L$ and $\bar a\in
W(K_{\omega})$, we have that $\ \Pi_{\omega}(\bar a)\in W(L)$.

\end{lemma}

Without loss of generality we can take $V$ to be $P^{n}(L)$ and consider the case $m=2$. The Segre embedding is defined by;\\

$P^{n}(L)\times P^{n}(L)\rightarrow P^{n(n+2)}(L)$\\

$((x_{0}:\ldots :x_{n}),(y_{0}:\ldots :y_{n}))\mapsto (x_{0}y_{0}:\ldots:x_{0}y_{n}:x_{1}y_{0}:\ldots :x_{n}y_{n})$\\

and the following diagram is easily checked to commute:\\

\begin{eqnarray*}
\begin{CD}
P^{n}(K_{i+1})\times P^{n}(K_{i+1})@>Segre>>P^{n(n+2)}(K_{i+1})\\
@VV\pi_{i+1}\times\pi_{i+1} V  @VV\pi_{i+1} V\\
P^{n}(K_{i})\times P^{n}(K_{i})@>Segre>>P^{n(n+2)}(K_{i})\\
\end{CD}
\end{eqnarray*}

Therefore, it is sufficient to prove that the property holds for $\pi_{i+1}:P^{n(n+2)}(K_{i+1})\rightarrow P^{n(n+2)}(K_{i})$. This is trivial to check using the fact that the residue map on $K_{i}[[t^{1/m}]]$ is a ring homomorphism fixing $K_{i}$.\\

Now we show that $(P^{n}(K_{\omega}),\Pi_{\omega})$ has the following universal property;\\

\begin{lemma}

Let $L\subset L_{m}$ be a field extension of transcendence degree
$m$, $V$ a smooth projective variety defined over $L$ and suppose
the map $\pi:V(L_{m})\rightarrow V(L)$ is given satisfying the
conclusion of Lemma 2.1. Then there exists an $L$-embedding
$\alpha_{L}:L_{m}\rightarrow K_{\omega}$ with the property that
$\Pi_{\omega}\circ \alpha_{L}=\pi$.

\end{lemma}

Choose a transcendence basis $\{t_{1},\ldots t_{m}\}$ for $L_{m}$ over $L$. We may assume that
$V$ is $P^{n}(L)$ for some $L$ and that $\pi:P^{n}(L_{m})\rightarrow P^{n}(L)$ is defined as above for some discrete valuation $v$ on $L_{m}$ with residue field $L$. Altering $(t_{1},\ldots t_{m})$ by automorphism if necessary, we may assume that $v|L(t_{1},\ldots,t_{m})$ is the canonical valuation given by;\\

$v_{res}:L(t_{1},\ldots,t_{m})\rightarrow ({\mathcal Z}^{m},<)$\\

$v_{res}(t_{1}^{i_{1}}\ldots t_{m}^{i_{m}})=(i_{1}\ldots i_{m})$\\

where $<$ denotes the lexographic ordering on ${\mathcal Z}^{m}$.
The completion of $L(t_{1},\ldots,t_{m})$ with respect to
$v_{res}$ is the formal Laurent series in $m$ variables
$(L((t_{1},\ldots,t_{m})),\overline{ v_{res}})$. Let $(\hat
{L_{m}},\bar v)$ be the completion of $L_{m}$ with respect to $v$,
then, as $L((t_{1},\ldots t_{m}))$ is Henselian with respect to
$\bar v_{res}$, $\bar v$ is the unique extension of $\bar v_{res}$
to $\hat {L_{m}}$. Now, for $1\leq i\leq m$, there exist canonical
isomorphisms  between $L((t_{i},\ldots,t_{m}))$ and
$L((t_{i}))((t_{i+1},\ldots,t_{m}))$. These combine to give an
$L$-embedding $\alpha_{L}$ of $L((t_{1},\ldots,t_{m}))$ into
$K_{m}$. Moreover, an easy calculation shows that
$\Pi_{\omega}\circ \alpha_{L}=\pi$ on $L(t_{1},\ldots,t_{m})$.
Now, by the uniqueness of the valuation extension from
$\alpha(L((t_{1},\ldots,t_{m})))$ to
$K_{m}=\alpha(L((t_{1},\ldots,t_{m})))^{alg}$, for any extension
of $\alpha$ to an embedding of $L_{m}$ into
$K_{m}$, we have $\Pi_{\omega}\circ \alpha_{L}=\pi$ on $L_{m}$ as well.\\

\end{section}

\begin{section}{Infintesimal Neighborhoods}

From now on, we fix a pair of Zariski structures and the specialisation map $\Pi_{\omega}$, to give a triple $((V(L),V(K_{\omega}),\Pi_{\omega})$ where $V$ is a smooth projective variety defined over $L$. \\

\begin{defn}
For  $\bar a\in{V(L)}^{n}$, we define the infintesimal
neighborhood of $\bar
a$ to be;\\

${\mathcal{V}}_{\bar a}=\Pi_{\omega}^{-1}(\bar a)$\\
\end{defn}

The first property of infintesimal neighborhoods is that we can
move inside closed sets.

\begin{lemma}

If $W(\bar y)$ is an irreducible closed set defined in $V(L)$,
$\bar b\in W$ and $dim(W)=r$, then there exists a $\bar
b'\in{\mathcal{V}}_{\bar b}\cap W(K_{\omega})$ such that $dim(\bar
b'/L)=r$

\end{lemma}

\begin{proof}

Consider the collection of constructible sets inside $V(L)^{n}$\\

$W(\bar y)\cup\{\neg C(\bar y): C\ closed,\ definable\ over\ $L$, dim (W(\bar y)\cap C(\bar y)<r)\}$\\

As $W$ is irreducible of dimension $r$, any finite subcollection
has a realisation in $V(L)^{n}$. By compactness, we can find a
realisation $\bar b'$ in $W(K)$ for $L\subset K$ such that
$dim(\bar b'/L)=r$. It then follows that we can define a partial
specialisation $\pi:V(K)\rightarrow V(L)$ by setting $\pi(\bar
b')=\bar b$, for if $C(\bar y)$ is a closed set defined over $L$
such that $\neg C(\bar b)$, then we must have that $dim (W(\bar
y)\cap C(\bar y)<d$ otherwise, $W$ being irreducible, $W(\bar
y)\subset C(\bar y)$, so by construction $\neg C(\bar b')$ also
holds. Now, using Lemma 2.2 applied to the field $L(\bar b')$
which has transcendence degree $r$ over $L$, we may assume that
$L(\bar b')\subset K_{r}\subset K_{\omega}$ and the specialisation
$\pi$ is given by the restriction of $\Pi_{\omega}$.
\\
\end{proof}

We now come to the critical theorem, a more general version of which was originally proved by Zilber in the context of abstract Zariski structures;\\

\begin{theorem}

Suppose that $F\subset D\times V^{k}$ is an irreducible finite
cover of $D$ with $D$ a smooth subvariety of $V^{m}$, and $F,D$
defined over $L$, such that $F(a,b)$. If
$a'\in{\mathcal{V}}_{a}\cap D(K_{\omega})$ is generic in $D$ over
$L$, then we can find $b'\in{\mathcal{V}}_{b}$ such that
$(a',b')\in F(K_{\omega})$.

\end{theorem}

We here only sketch the proof, full details may be found in \cite{Pez}. We first consider the following collection of constructible sets defined over $K_{\omega}$, with $a'\in{\mathcal{V}}_{a}\cap D(K_{\omega})$ generic over $L$;\\

$\{F(a',y)\}\cup\{\neg C(d,y): d\in V(K_{\omega}), \neg C(\Pi_{\omega}(d),b)\}$\\

As $F$ is a finite cover and $K_{\omega}$ is algebraically closed, a realisation $b'$ of this collection lies in $V(K_{\omega})$ and $F(a',b')$ holds. Moreover, $\Pi_{\omega}(b')=b$, otherwise, as the diagonal $x=y$ is closed, we have that $b'\neq y$ is in the collection which is ridiculous. \\

If the collection is inconsistent, we find a closed set $Q\subset V^{n+k}$ such that $F(a',y)\subseteq Q(d,y)$ whereas $\neg Q(\pi(d),b)$. \\

The point of the smoothness assumption is to show that the parameter space \\

$L(x,z)\subset D\times V^{n}\ =\{(x,z):F(x,y)\subset Q(z,y)\}$\\

which in general is not relatively closed in $D\times M^{n}$ at least corresponds  to a closed set over a dense open subset of $D$. More precisely, there is a closed subvariety $P(x,z)\subset D\times M^{n}$ and $D'\subset D$, $dim(D')<dim(D)$, all defined over $L$, such that \\

1. $P(x,z)\subset L(x,z)$.\\

2. $L(x,z)\subset P(x,z)\cup (D'\times V^{n})$  $(*)$\\

We have by assumption that $L(a',d)$ holds. As $a'$ was chosen to be generic over $L$ and $a'\notin D'$, $P(a', d)$ holds. Applying the specialisation $\Pi_{\omega}$ gives that $P(a,\Pi_{\omega}(d))$, hence $F(a,y)\subset Q(\Pi_{\omega}(d),y)$, hence $Q(\Pi_{\omega}(d),b)$ holds as well, contradicting the assumption. \\

\begin{rmk}

In fact the theorem can be improved to give the following more
general result;

Suppose that $F\subset D\times V^{k}$ is an irreducible
generically finite cover of $D$ with $D$ a subvariety of $V^{m}$.
Then, if $a\in D$ is a regular point for the cover and contained
in the non-singular locus of $D$, $a'\in{\mathcal V}_{a}\cap
D(K_{\omega})$ is generic in $D$ over $L$, then we can find
$b'\in{\mathcal V}_{b}$ such that $(a',b')\in F(K_{\omega})$.

\end{rmk}

\end{section}

\begin{section}{Zariski Unramified Maps and Multiplicity}

The purpose of introducing infintesimal neighborhoods is to define an abstract notion of Zariski
multiplicity. \\

\begin{defn}{Zariski multiplicity}\\

Let hypotheses be as in Theorem 3.3\\

Given $(a,b)\in F$, set\\

   $mult_{ab}(F/D)=Card (F(a',K_{\omega}))\cap{\mathcal{V}}_{b})\ for\ a'\in{\mathcal{V}}_{a}\cap D\ generic\ over\ L$\\
\end{defn}

We want to show this is well defined.\\
\begin{proof}

Suppose $a''\in{\mathcal{V}}_{a}\cap D$ with $Card(F(a'',K_{\omega})\cap{\mathcal{V}}_{b})=n$. Consider the relation $N(x,y_{1},\ldots,y_{n})\subset D\times V^{nk}$, given by\\

$N(x,y_{1},\ldots,y_{n})=F(x,y_{1})\wedge\ldots\wedge F(x,y_{n})$\\

Then we have that $N$ is a finite cover of $D$ and, moreover, by smoothness of $D$, each irreducible component of $N$ has dimension at least\\

$n(dim(F)+(n-1)k)-(n-1)(dim(D)+nk)=dim(D)+n(n-1)k-n(n-1)k=dim(D)$\\

so clearly each component is a finite cover of $D$. Now, choose an irreducible component $N_{i}$ containing $(a'',b''_{1},\ldots,b''_{n})$, so by specialisation also contains $(a,b,\ldots,b)$ and consider the open set $U\subset N_{i}$ given by\\

$U(x,y_{1},\ldots,y_{n})=N_{i}(x,y_{1},\ldots,y_{n})\wedge y_{1}\neq y_{2}\neq\ldots\neq y_{n}$\\

Then, for $a'\in{\mathcal{V}}_{a}$ generic in $D$, it follows we can find a tuple $(b'_{1},\ldots, b'_{n})$ such that $N_{i}(a',b'_{1},\ldots,b'_{n})$, and $(b'_{1},\ldots,b'_{n})\in{\mathcal{V}}_{(b,\ldots,b)}$. As is easily checked, the tuple $(a',b'_{1},\ldots,b'_{n})$ is generic inside $N_{i}$, hence must lie inside $U$. This proves that the $b'_{1},\ldots,b'_{n}$ are distinct, hence $Card(F(a',K_{\omega})\cap {\mathcal V}_{b})\geq n$.\\
\end{proof}

\begin{defn}
We say that a point $(ab)\in F$ is Zariski ramified if $mult_{ab}(F/D)\geq 2$. Otherwise, we call such a point Zariski unramified. \\
\end{defn}

Now suppose $F\subset D\times V^{n}$ is an irreducible  finite cover of $D$ with $D$ smooth, then we have the following easily checked lemma\\

\begin{lemma}

$mult_{a}(F/D)=_{def}\Sigma_{b\in F(a,L)}mult_{ab}(F/D)$ does not depend on the choice of $a\in D$, and  is equal to the size of a generic fibre over $D$\\

\end{lemma}

A simple consequence is the following:\\

\begin{lemma}

If $ \bar a'\in D(L)$, then $F(\bar a')$ contains a point of
ramification in the sense of Zariski structures iff $|F(\bar
a')|<|F( \bar a)|$ where $ \bar a$ is generic in $D$.

\end{lemma}

\begin{proof}
We have seen that $|F(\bar a)|=\Sigma_{\bar b\in F(\bar a',L)}mult_{\bar a',\bar b}(F/D)$.  If $|F(\bar a')|<|F(\bar a)|$, then there must exist $\bar b\in F(\bar a')$ with $mult_{(\bar a',\bar b)}(F/D)\geq 2$ so the result follows by the definition of ramification in Zariski structures. The converse is similar.\\
\end{proof}

We will also require the following results, that Zariski
multiplicity is multiplicative over composition and preserved by
open maps.

\begin{lemma}

Suppose that $F_{1},F_{2}$ and $F_{3}$ are smooth, irreducible,
with $F_{2}\subset F_{1}\times V^{k}$ and $F_{3}\subset
F_{2}\times V^{l}$ finite covers. Let $(abc)\in F_{3}\subset
F_{1}\times V^{k}\times V^{l}$. Then
$mult_{abc}(F_{3}/F_{1})=mult_{ab}(F_{2}/F_{1})mult_{bc}(F_{3}/F_{2})$.

\end{lemma}

\begin{proof}
To see this, let $m=mult_{ab}(F_{2}/F_{1})$ and $n=mult_{bc}(F_{3}/F_{2})$. Choose $a'\in {\mathcal V}_{a}\cap F_{1}(K_{\omega})$ generic over $L$. By definition, we can find distinct $b_{1}\ldots b_{m}$ in $V^{k}(K_{\omega})\cap{\mathcal V_{b}}$ such that $F_{2}(a',b_{i})$ holds. As $F_{2}$ is a finite cover of $F_{1}$, we have that $dim(a'b_{i}/L)=dim(a'/L)=dim(F_{1})=dim(F_{2})$, so each $(a'b_{i})\in {\mathcal V}_{ab}\cap F_{2}$ is generic over $L$. Again by definition, we can find distinct $c_{i1}\ldots c_{in}$ in $V^{l}(K_{\omega})\cap{\mathcal V_{c}}$ such that $F_{3}(a'b_{i}c_{ij})$ holds. Then the $mn$ distinct elements $(a'b_{i}c_{ij})$ are in ${\mathcal V_{abc}}$, so by definition of multiplicity $mult_{abc}(F_{3}/F_{1})=mn$ as required.\\
\end{proof}

\begin{lemma}

Let $\pi_{1}:X\rightarrow D$ and $\pi_{2}:Y\rightarrow D$ be
covers with assumptions as in remarks following Theorem 3.3.
Suppose moreover that there exist open smooth subvarieties
$U\subset X$ and $V\subset Y$ and an isomorphism $f:U\rightarrow
V$ such that $\pi_{2}\circ f=\pi_{1}$ on $U$. Then if $a\in D$ is
a regular point for the cover $\pi_{1}$ and $(ab)\in U$,
$mult_{ab}(X/D)=mult_{af(b)}(Y/D)$.

\end{lemma}

\begin{proof}

We may assume that the open set $U$ is maximal with the property
that $(ab)\in U$ and there exists an isomorphism with $V\subset
Y$. Suppose $mult_{ab}(X/D)=m$. Then we can find $a'\in {\mathcal
V}_{a}\cap D(K_{\omega})$ generic in $D$ over $L$ and
$b_{1},\ldots,b_{m}$ distinct such that $X(a'b_{i})$ holds for
$1\leq i\leq m$. It will be sufficient to show that
$Y(a'f(b_{i}))$ holds and $f(b_{i})\in{\mathcal V}_{f(b)}$, for
$1\leq i\leq m$, then, as $f$ is injective,
$mult_{a,f(b)}(Y/D)\geq m$ and the result follows by symmetry. By
the fact that $\pi_{2}\circ f=\pi_{1}$ on $U$ we clearly have that
$Y(a'f(b_{i}))$ holds. Let $\overline{graph(f)}$ be the projective
closure of the graph of $f$ in the projective variety $X\times Y$
and $\pi_{X},\pi_{Y}$ the projections onto the coordinates $X$ and
$Y$. Then $\pi_{X}$ satisfies the conditions of the remarks after
Theorem 3.3, and moreover by assumption the point $(ab)\in X$ is
regular for the cover $\pi$ and contained in the non-singular
locus of $X$. Hence, we can find $(cd)\in{\mathcal V}_{af(b)}$
such that $\overline{graph(f)}(a'b_{i},cd)$ holds. As $graph(f)$
is a $1-1$ correspondence between $U$ and $V$, if
$(a'b_{i},cd)\in\overline{graph(f)}\setminus graph(f)$ then
$(a'b_{i},cd)\in F_{X}\cup F_{Y}$ where $F_{X}, F_{Y}$ consist of
the infinite fibres of the projections $\pi_{X}$ and $\pi_{Y}$
respectively. By $(DF)$, both of these are defined over $L$ and
have dimension strictly less than $graph(f)$. This contradicts the
fact that $(a'b_{i},cd)$ is generic inside $\overline {graph(f)}$
over $L$, hence $(a'b_{i},cd)\in graph(f)$ and as $f$ is a
bijection $(cd)=(a'f(b_{i}))$. This shows that
$f(b_{i})\in{\mathcal V}_{f(b)}$ as required.

\end{proof}

\end{section}

\begin{section}{Etale Morphisms and Algebraic Multiplicity}

We review here the algebraic notions which will be required in the following section. \\

\begin{defn}
A morphism $f$ of finite type between varieties $X$ and $Y$ is said to be etale if for all $x\in X$ there are open affine neighborhoods $U$ of $x$ and $V$ of $f(x)$ with $f(V)\subset U$ such that restricted to these neighborhoods the pull back on functions is given by the inclusion;\\

$f^{*}:L[V]\rightarrow L[V]{[x_{1},\ldots, x_{n}]\over f_{1},\ldots, f_{n}}$\\

and $det({\partial f_{i}\over \partial x_{j}})(x)\neq 0\ ,(*)$\\

\end{defn}

The coordinate free definition of etale is that $f$ should be flat and unramified.\\

The notion of an etale morphism simplifies considerably when we
assume that $X$ and $Y$ are smooth algebraic varieties over $L$,
see \cite{Mum};

\begin{theorem}

If $X$ and $Y$ are non-singular algebraic varieties over $L$ and
$f:X\rightarrow Y$ is a morphism, then $f$ is etale iff
$df:(m_{x}/m_{x}^{2})^{*}\rightarrow (m_{f(x)}/m_{f(x)}^{2})^{*}$
is an isomorphism everywhere.

\end{theorem}

\begin{rmk}

This gives us a convenient test for etaleness given an arbitrary morphism of finite type between smooth varieties $X$ and $Y$.  If we take local uniformisers $g_{1},\ldots g_{n}$ at $x\in X$, the $dg_{i}$ generate $\Omega_{X}$ freely on an open $U'$ of $x$.  If we pull back a set of uniformisers $f^{*}f_{1},\ldots,f^{*}f_{n}$ on $Y$ to $X$, we can locally define the Jacobian $Jac^{\bar f}_{\bar g}$ to be;\\

$det({\partial f^{*}f_{i}\over\partial g_{j}})$\\

which means write the $1$-forms
$f^{*}df_{i}=\Sigma_{j}a_{ij}dg_{j}$ and take $det(a_{ij})$. If
$f$ is etale in a neighborhood of $x$, the $f^{*}df_{i}$ also
generate $\Omega_{X}$ freely on an open $U''$ of $x$. Taking the
intersection $U''=U\cap U'$, gives us that the Jacobian $Jac^{\bar
f}_{\bar g}|U''\neq 0$. Conversely, if $Jac^{\bar f}_{\bar
g}(x)\neq 0$, then it is non zero on an open neighborhood $U''$ of
$x$ and by the above theorem we have that $f$ is etale
 on this neighborhood.
\end{rmk}

We will also require some facts about the etale topology on an
algebraic variety $Y$. We consider a category $Y_{et}$ whose
objects are etale morphisms $U\rightarrow Y$ and whose arrows are
$Y$-morphisms from $U\rightarrow V$. This category has the
following $2$ desirable properties. First given $y\in Y$, the set
of objects of the form $(U,x)\rightarrow (Y,y)$ form a directed
system, namely $(U,x)\subset (U',x')$ if there exists a morphism
$U\rightarrow U'$ taking $x$ to $x'$. Secondly, we can take
``intersections'' of open sets $U_{i}$ and $U_{j}$ by considering
$U_{ij}=U_{i}\times_{Y}U_{j}$; the projection maps are easily show
to be etale and the composition of etale maps is etale, so
$U_{ij}\rightarrow Y$ still lies in $Y_{et}$. (Note that we can
develop the theory of etale cohomology for an arbitrary Zariski
structure, this will be a subject of further investigation) If $Y$
is an irreducible variety over $K$, then all etale morphisms into
$Y$ must come from reduced schemes of finite type over $K$, though
they may well fail to be irreducible considered as algebraic
varieties. Now we can define the local ring of $Y$ in the etale toplogy to be;\\

$O_{y,Y}^{\wedge}=lim_{\rightarrow, y\in U}O_{U}(U)$\\

As any open set $U$ of $Y$ clearly induces an etale morphism $U\rightarrow_{i}Y$ of inclusion, we have that $O_{y,Y}\subset O_{y,Y}^{\wedge}$. We want to prove that $O_{y,Y}^{\wedge}$ is a Henselian ring and in fact the smallest Henselian ring containing $O_{y,Y}$. We need the following lemma about Henselian rings;\\

\begin{lemma}

Let $R$ be a local ring with residue field $k$. Suppose that $R$ satisfies the following condition;\\

If $f_{1},\ldots f_{n}\in R[x_{1},\ldots x_{n}]$ and $\bar f_{1}\ldots \bar f_{n}$ have a common root $\bar a$ in $k^{n}$, for which $Jac(\bar f)(\bar a)=({\partial \bar f_{i}\over \partial x_{j}})_{ij}(\bar a)\neq 0$, then $\bar a$ lifts to a common root in $R^{n}$\ (*).\\

Then $R$ is Henselian.\\

\end{lemma}

It remains to show that $O_{y,Y}^{\wedge}$ satisfies $(*)$.

\begin{proof}

Given $f_{1,}\ldots f_{n}$ satisfying the condition of $(*)$, we
can assume the coefficients of the $f_{i}$ belong to
$O_{U_{i}}(U_{i})$ for covers $U_{i}\rightarrow Y$; taking the
intersection $U_{1\ldots i\ldots n}$ we may even assume the
coefficients define functions on a single etale cover $U$ of $Y$.
By the remarks above we can consider $U$ as an algebraic variety
over $K$, and even an affine algebraic variety after taking the
corresponding inclusion. We then consider the variety $V\subset
U\times A^{n}$ defined by $Spec({R(U)[x_{1},\ldots,x_{n}]\over
f_{1},\ldots f_{n}})$. Letting $u\in U$ denote the point in $U$
lying over $y\in Y$, the residue of the coefficients of the
$f_{i}$ at $u$ corresponds to the residue in the local ring $R$,
which tells us exactly that the point $(u,\bar a)$ lies in $V$. By
the Jacobian condition, we have that the projection
$\pi:V\rightarrow U$ is etale at the point $(u,\bar a)$, and hence
on some open neighborhood of $(u,\bar a)$, using
Nakayama's Lemma applied to $\Omega_{V/U}$. Therefore, replacing $V$ by the open subset
 $U'\subset V$ gives an etale cover of $U$ and therefore of $Y$, lying over $y$. Now clearly
 the coordinate functions $x_{1},\ldots x_{n}$ restricted to $U'$ lie in
  $O_{y,Y}^{\wedge}$ and lift the root $\bar a$ to a root in $O_{y,Y}^{\wedge}$\\
\end{proof}

We define the Henselization of a local ring $R$ to be the smallest Henselian ring $R'\supset R$, with $R'\subset Frac(R)^{alg}$. We have in fact that;\\

\begin{theorem}

Given an algebraic variety $Y$, $O_{y,Y}^{\wedge}$ is the
Henselization of $O_{y,Y}$

\end{theorem}

\begin{defn}

Given smooth projective curves $C_{1}$, $C_{2}$ and a finite
morphism $f:C_{1}\rightarrow C_{2}$, the algebraic multiplicity of
$f$ at $a$ is $ord_{a}(f^{*}h)$ where $h$ is a local uniformiser
for $C_{2}$ at $f(a)$.

\end{defn}

\begin{rmk}
This is independent of the choice of $h$, as the quotient of $2$ uniformisers $h/h'$ is a unit in ${\mathcal{O}}_{f(a)}$. Given finite morphisms $f:C_{3}\rightarrow C_{2}$ and $g:C_{2}\rightarrow C_{1}$, if $ord_{a,f(a)}(C_{3}/C_{2})=m$ and $ord_{f(a),gf(a)}(C_{2}/C_{1})=n$, then taking a local uniformiser $h$ at $gf(a)$, we have that $g^{*}h=h_{1}^{n}u$ locally at $f(a)$ for a unit $u$ and uniformiser $h_{1}$ in ${\mathcal O}_{f(a)}$. Similarily $f^{*}g^{*}h=h_{2}^{mn}u'$ for a unit $u'$ and  uniformiser$ h_{2}$ in ${\mathcal O}_{a}$. This shows that $ord_{a,gf(a)}(C_{3}/C_{1})=mn$, so the branching number is also multiplicative for smooth projective curves.  \\

\end{rmk}

\begin{defn}

Given smooth projective varieties $X_{1}$, $X_{2}$ and a finite
morphism $f:X_{1}\rightarrow X_{2}$, the algebraic multiplicity
$mult_{af(a)}^{alg}(X_{1}/X_{2})$ of $f$ at $a\in X_{1}$ is
$length(O_{a,X_{1}}/f^{*}m_{f(a)})$ where $m_{f(a)}$ is the
maximal ideal of the local ring $O_{f(a)}$.

\end{defn}

\begin{rmk}

Note that this is finite, by the fact that finite morphisms have
finite fibres and the ring $O_{a,X_{1}}/f^{*}m_{f(a)}$ is a
localisation of the fibre $f^{-1}(f(a))\cong
R(f^{-1}(U))\otimes_{R(U)}L\cong R(f^{-1}(U))/m_{f(a)}$ where $U$
is an affine subset of $X_{2}$ containing $f(a)$.

\end{rmk}

We now have the following, which generalises the result for
curves;\\

\begin{theorem}{Algebraic multiplicity is multiplicative};\\

Given finite morphisms $f:X_{3}\rightarrow X_{2}$ and
$g:X_{2}\rightarrow X_{1}$ between smooth projective varieties,
for $a\in X_{3}$ we have that\\

$mult_{af(a)}(X_{3}/X_{2})mult_{f(a)gf(a)}(X_{2}/X_{1})=mult_{agf(a)}(X_{3}/X_{1})$.\\

\end{theorem}

\begin{proof}

The proof is an exercise in algebra, which we give for want of a
convenient reference. First, the morphisms $f$ and $g$ are flat.
This requires the following lemma, given as an exercise in
\cite{Hart}, and the fact that smooth varieties are regular and Cohen-Macauley;\\

\begin{lemma}

Let $f:X\rightarrow Y$ be a morphism of varieties over $L$. Assume
that $Y$ is regular, $X$ is Cohen-Macauley and that every fibre of
$f$ has dimension equal to $dim(X)-dim(Y)$. Then $f$ is flat.

\end{lemma}

Now we have a tower of local rings $(R,m)\subset (S,n)\subset
(T,o)$ with algebraically closed residue field $L$. Each extension
is free by the flatness result and finiteness. For a finite free
extension $(R,m)\subset (S,n)$ of local rings, we also have the easily checked result that;\\

$[S:R]=dim_{Fr(R)}S\otimes_{R}
Fr(R)=dim_{L}S\otimes_{R/mR}R/mR=dim_{L}(S/mS) (*)$.\\

 For an extension $(R,m)\subset (S,n)$ of local rings, we have that
$length(S/mS)=dim_{L}(S/mS)$, hence, by (*), the theorem reduces
to checking that $[T:R]=[T:S][S:R]$ which is standard.\\

\end{proof}

\end{section}

\begin{section}{Equivalence of the Notions}

This section is devoted to the main proofs of the paper, namely that the notions developed in Sections 4 and 5 are essentially equivalent for morphisms between smooth projective varieties. \\

\begin{theorem}

Let hypotheses be as in Theorem 3.3, with the additional
assumption that $char L=0$, then $F$ is a Zariski unramified cover
of $D$ iff $F$ is an etale cover of $D$.

\end{theorem}

Let $pr$ be the projection map of $F$ onto $D$, then $pr$ is a projective morphism. By Zariski's Main Theorem, $pr$ factors as a composition $F\rightarrow_{pr_{1}} F'\rightarrow_{pr_{2}} D$ with $pr_{1}$ having connected fibres, $pr_{1*}F=F'$and $pr_{2}$ a finite morphism. The formal inverse $pr_{1}^{-1}$ from $F'$ to $F$ is a morphism corresponding to the identification of $pr_{1*}F$ and $F'$, hence $pr_{1}$ is in fact an isomorphism. We may therefore assume that $pr$ is a finite morphism.\\

Now suppose that $pr$ is etale, then, $pr$ is flat, see \cite{Mum} for how this follows from Definition 5.1. As $D$ is irreducible, \\

$dim_{k(y)}(f_{*}(O_{F})\otimes_{O_{y}}k(y))$\\

is independent of $y\in D$. As $pr$ is etale, $pr_{*}:T_{x,F}\rightarrow T_{pr(x),D}$ is an isomorphism, hence, by a simple calculation;\\

$dim_{k(y)}(f_{*}(O_{F})\otimes_{O_{y}}k(y))=|F(y)|$ for $y\in D$.\\

This shows that $|F(y)|$ is independent of $y\in D$. By Lemma 4.4, this shows that $pr$ is a Zariski unramified cover.\\

Conversely, suppose that $pr$ is Zariski unramified. We first show
that for generic $\bar a\in D$, $|F(\bar
a)|=deg(pr)=deg[k(F):k(D)]$. As $char(k(F))=0$, the extension is
seperable so we can find a primitive element $g\in k(F)$ such that
$k(F)=k(D)(g)$. Clearly the minimum polynomial $p$ of $g$ over
$k(D)$ has degree $n=deg[k(F):k(D)]$. Let $h_{1},\ldots h_{n-1}\in
k(D)$ be the coefficients of $p$, then $R(D)(h_{1}\ldots h_{n-1})$
determines the function ring of a Zariski open subset $U$ of $D$.
Clearly $R(U)[g]$ is an integral extension of $R(U)$ and
corresponds to the projection restricted to $U'=pr^{-1}(U)\cap
g\neq 0$. By dimension theory, the zero set $Z(g)\subset D$ cannot
intersect with a generic fibre of the original map
$pr:F\rightarrow D$. Now we consider the discriminant $D(p)$ of
the polynomial $p$ as a regular function on $U$ and we have that
for generic $\bar a\in U$ that $D(p)(\bar a)\neq 0$. This implies
that for generic $\bar a\in U
$ $|pr^{-1}(\bar a)|=n=deg[k(F):k(D)]$. Now we are in a position to apply Theorem $5$, p145, of
 \cite{Shaf} which requires that $D$ should be smooth, namely that $pr_{*}:T_{x,F}\rightarrow T_{pr(x),D}$ is an isomorphism for $x\in F$. As $F$ and $D$ were assumed to be nonsingular, this is sufficient to show that $pr$ is etale by Theorem 5.2.\\

\begin{rmk}
                                                  When $char(L)=p$, the analogy fails.  If we consider the Frobenius map $Fr:P^{1}\rightarrow P^{1}$, then $Graph(Fr)\subset P^{1}\times P^{1}$ is a finite cover of $P^{1}$ and both $Graph(Fr)$ and $P^{1}$ are smooth. The projection map $pr$ onto the second coordinate is unramified in the sense of Zariski structures as $pr$ is a bijection. However $pr$ fails to be etale in the sense of algebraic geometry as $pr_{*}:T_{x,Graph(Fr)}\rightarrow T_{pr(x),P^{1}}$ is zero everywhere. However the following theorem shows that this is the only bad example and highlights one advantage of the Zariski method, namely that it is insensitive \emph{only} to Frobenius.

\end{rmk}

\begin{theorem}

Let hypotheses be as in Theorem 3.3, with the additional
assumption that $char(L)=p\neq 0$. If $F$ is an etale cover of
$D$, $F$ is a Zariski unramified cover. Conversely, if $F$ is a
Zariski unramified cover, then $pr$ factors as a composition
$F\rightarrow_{pr_{1}} F'\rightarrow_{pr_{2}} D$ in $Proj$ with
$pr_{1}$ a purely inseperable connected cover and $pr_{2}$ an
etale cover.

\end{theorem}

\begin{proof}

As in the previous theorem, we may assume that $pr$ is a finite
morphism.
Suppose first that $F\rightarrow D$ is a finite morphism with $F$ and $D$ affine.  We  first find a field $L$ such that $k(F)/L$ is a purely inseperable extension and $L/k(D)$ is seperable. Let $R'$ be the integral closure of $R(D)$ in $L$ and $R''$ the integral closure of $R(D)$ in $k(F)$. As $R(F)$ is integral over $R(D)$ we have that $R(F)\subset R''$, but $F$ was assumed to be smooth so $R(F)$ is integrally closed in $k(F)$ and therefore $R''=R(F)$. As the extensions $k(D)\subset L\subset k(F)$ are finite algebraic, by \cite{ZarSam}, both $R(F)$ and $R'$ are finite $R'$ and $R(D)$ modules respectively. Therefore, corresponding to the ring inclusions \\

$R(D)\rightarrow R'\rightarrow R(F)$\\

we have the sequence of finite morphisms\\

$F\rightarrow_{pr_{1}} Spec(R')\rightarrow_{pr_{2}} D$\\

We first consider the cover $F\rightarrow_{pr_{1}} Spec(R')$. Let
$g_{1},\ldots g_{m}$ generate $R(F)$ over $R'$. As the extension
$k(F)/L$ is purely inseperable, we can write the minimum
polynomials $p_{i}$ of $g_{i}$ in the form
$r_{i,0}g^{p^{n_{i}}}-r_{i,1}=0$ where $r_{i,0}$ and $r_{i,1}$ are
in $R'$. As $R(F)/R'$ is finite, we can also find monic
polynomials $q_{i}$ with coefficients in $R'$ satisfied by
$g_{i}$. Choose polynomials
$t_{i}=s_{i,0}x^{m_{i}}+s_{i,1}x_{m_{i}-1}+\ldots s_{i,m_{i}}$
such that $p_{i}t_{i}=q_{i}$. By equating coefficients, we have
that $r_{i,0}=s_{i,0}^{-1}$ and $r_{i,1}/r_{i,0}\in R'$. Hence, we
can take the $p_{i}$ to be monic with coefficients in $R'$. As the
$p_{i}$ are minimal monic polynomials, we conclude that that
$R(F)$ is an extension of the form
$R'[g_{1},\ldots,g_{m}]/(g_{1}^{p^{n_{1}}}-\lambda_{1},\ldots,
g_{m}^{p^{n_{m}}}-\lambda_{m})$ with $\lambda_{i}\in R'$. This is
easily checked to be a connected cover of $Spec(R')$. In fact if
we let $\theta=(Fr^{-n_{1}},\ldots, Fr^{-n_{m}})\circ
(\lambda_{1}\ldots \lambda_{m})$, where the $\lambda_{i}$ are
considered as regular functions on $Spec(R')$ and $Fr^{-n_{i}}$ is
the formal inverse Frobenius map, then the cover corresponds to
the projection of $Graph(\theta)\subset Spec(R')\times A^{m}$ onto
$Spec(R')$. As $F$ was assumed to be smooth, $Spec(R')$ is a
smooth seperable Zariski unramified cover of $D$. Applying the
previous theorem, we conclude that $Spec(R')$ is an etale cover of
$D$. Now, for the case when $F$ and $D$ are projective varieties,
let $U_{i}$ be an affine cover of $D$ and $R'(U_{i})$ the
corresponding normalisations. By uniqueness of integral closure,
the $R'(U_{i})$ patch to form a cover $F'$ of $D$. In fact, by a
classical result, see \cite{Mum}, we may assume that $F'$ is a
smooth projective variety. As etaleness is a local condition for
smooth varieties, the cover $F'$ is etale. Finally, check that the
local maps $pr_{1}:F_{i}\rightarrow R'(U_{i})$ patch on overlaps
to give a morphism $pr_{1}:F\rightarrow F'$. Clearly, this is an
insperable connected cover, in fact if $F'$ is defined by the
homogenous equations $<f_{1},\ldots f_{n}>$ inside $P^{N}$, then
$F$ is isomorphic to the closed subvariety of $P^{N}\times P^{m}$
defined by the extra equations
$<Y_{i}^{p^{n_{i}}}X_{N}^{j(i)}-\lambda_{i}(X_{0},\ldots,X_{N})Y_{0}^{p^{n_{i}}}>$
where $1\leq i\leq m$ and $j(i)$ is the degree of the polynomial
$\lambda_{i}$ in the affine coordinates $P^{N}(L)_{i}$.

\end{proof}

\begin{rmk}
We now show that the notions of Zariski multiplicity and algebraic multiplicity coincide when $char(L)=0$, as usual with assumptions being as in Theorem 3.3., and find an anlogous result when $char(L)=p$.  Unfortunately, it does not seem possible to achieve this by counting points in the fibres, as in the previous theorems, so we need to find a local method. This will be the subject of the remainder of this section.\\

\end{rmk}

For ease of exposition, we first consider the case when $F$ and $D$ are curves. We will point out the necessary modifications  for the case when $F$ and $D$ are arbitrary smooth projective varieties in the next theorem.\\

\begin{theorem}

Let hypotheses be as in Theorem 3.3, with the additional
assumption that $char(L)=0$ and $F$, $D$ are curves. Then the
notions of Zariski multiplicity and algebraic multiplicity
coincide.

\end{theorem}

\begin{proof}

As $D$ has a non-constant meromorphic function, we can write $D$ as a finite cover of $P^{1}(L)$.  As we have checked both algebraic multiplicity and Zariski multiplicity are multiplicative over composition, a straightforward calculation shows that we need only check the notions agree for the branched finite cover $\pi:F\rightarrow P^{1}(L)$. (1)\\

Now consider this cover restricted to $A^{1}$, let $x$ be the canonical cooordinate  with $ord_{a} (\pi^{*}(x))=m$, so we have that $\pi^{*}x=h^{m}u$ , for $u$ a unit in ${\mathcal{O}}_{a}$ and $h$ a uniformiser at $a$. (2)\\

As $u$ is a unit and $char(L)=0$, the equation $z^{m}=u$ splits in the residue field of ${\mathcal O}^{\wedge}_{a}$. By Hensel's Lemma and Theorem 5.5, it is solvable in ${\mathcal O}_{a}^{\wedge}$. By the definition of ${\mathcal{O}}_{a}^{\wedge}$, we can find an etale morphism $\pi:(U,b)\rightarrow (F,a)$ containing such a solution in the local ring ${\mathcal{O}}_{b}$. We may assume that $U$ is irreducible and moreover, as $\pi$ is etale, that $U$ is smooth. (3)\\

Now we can embed $U$ in a projective smooth curve $F'$ and, as $F$ is smooth,  extend the morphism $\pi$ to a projective morphism from $F'$ to $F$. (4)\\

We claim that $(ba)\in graph(\pi)\subset F'\times F $ is unramified in the sense of Zariski structures. For this we need the following fact whose algebraic proof relies on the fact that etale morphisms are flat, see \cite{Milne};\\

\begin{fact}

Any etale morphism can be locally presented  in the form \\

\begin{eqnarray*}
\begin{CD}
V@>g>>Spec((A[T]/f(T))_{d})\\
@VV\pi V  @VV\pi' V\\
U@>h>>Spec(A)\\
\end{CD}
\end{eqnarray*}

where $f(T)$ is a monic polynomial in $A[T]$, $f'(T)$ is invertible in $(A[T]/f(T))_{d}$ and $g,h$ are isomorphisms.   (5)\\
\end{fact}

Using Lemma 4.6 and the fact that the open set $V$ is smooth, we may safely replace $graph(\pi)$ by $\overline {graph (\pi')}\subset F''\times F$ where $F''$ is the projective closure of $Spec((A[T]/f(T))$, $F$ is the projective closure of $Spec(A)$ and $\overline {graph(\pi')}$ is the projective closure of $graph(\pi')$ and show that $(g(b)a)$ is Zariski unramified. Note that over the open subset $U=Spec(A)\subset F$, $\overline{graph(\pi')}=Spec((A[T]/f(T)$ as this is closed in $U\times F''$.  For ease of notation, we replace $(g(b)a)$ by $(ba)$. (6)\\

Suppose that $f$ has degree $n$. Let $\sigma_{1}\ldots \sigma_{n}$ be the elementary symmetric functions in $n$ variables $T_{1},\ldots T_{n}$. Consider the equations\\

$\sigma_{1}(T_{1},\ldots, T_{n})=a_{1}$\\

$\ldots$\\

$\sigma_{n}(T_{1},\ldots,T_{n})=a_{n}$ (*)\\

where $a_{1},\ldots a_{n}$ are the coefficients of $f$ with
appropriate sign. These cut out a closed subscheme $C\subset
Spec(A[T_{1}\ldots T_{N}])$. Suppose $(ba)\in
graph(\pi')=Spec(A[T]/f(T))$ is ramified in the sense of Zariski
structures, then I can find $(a'b_{1}b_{2})\in {\mathcal V}_{abb}$
with $(a'b_{1})$,$(a'b_{2})\in Spec(A(T)/f(T))$ and $b_{1},b_{2}$
distinct. Then complete $(b_{1}b_{2})$ to an $n$-tuple
$(b_{1}b_{2}c_{1}'\ldots c_{n-2}')$ corresponding to the roots of
$f$ over $a'$. The tuple $(a'b_{1}b_{2}c_{1}'\ldots c_{n-2}')$
satisfies $C$, hence so does the specialisation $(abbc_{1}\ldots
c_{n-2})$. Then the tuple $(bbc_{1}\ldots c_{n-2})$ satisfies
$(*)$ with the coefficients evaluated at $a$. However such a
solution is unique up to permutation and corresponds to the roots
of $f$ over $a$. This shows that $f$ has a double root at $(ab)$
and therefore $f'(T)|_{ab}=0$. As $(ab)$ lies inside
$Spec(A[T]/f(T))_{d}$, this contradicts the fact that $f'$ is
invertible in $A[T]/f(T))_{d}$. (7)\\

In $(2)$ we may therefore assume that $\pi^{*}x=h^{m}$ for $h$ a local uniformiser at $a$. Now we have the sequence of ring inclusions given by \\

$L[x]\rightarrow L[x,y]/(y^{m}-x)\rightarrow R$\\

\ \ \ \ \ \ \ \ \ $x\mapsto \pi^{*}x, y\mapsto h$\\

where $R$ is the coordinate ring of $F$ in some affine neighborhood of $a$. It follows that we can factor our original map such that $F$ is etale near $a$ over the projective closure of $y^{m}-x=0$. (8)\\

Again, repeating the argument from (4) to (7), we just need to
check that the projective closure of $y^{m}-x$ has multiplicity
$m$ at $0$ considered as a cover of $P^{1}(\bar k)$. This is
trival, let $\epsilon\in {\mathcal V}_{0}$ be generic over
$\mathcal M$,then as we are working in characteristic $0$ we can
find distinct $\epsilon_{1},\ldots \epsilon_{m}$ in ${\mathcal
M}_{*}$ solving $y^{m}=\epsilon$. By specialisation, each
$\epsilon_{i}\in{\mathcal V}_{0}$. (9)
\end{proof}

\begin{theorem}

Let hypotheses be as in Theorem 6.5, with the modification that
$char(L)=p\neq 0$. If $e$ denotes the Zariski multiplicity and $d$
the algebraic multiplicity at $a\in F$, then $d=ep^{n}$ and $\pi$
factors as $F\rightarrow_{h}F'\rightarrow_{g}D$ with $h=Frob^{n}$
and $g$ having algebraic multiplicity $e$ at $h(a)$.

\end{theorem}

By Theorem 6.3, we can factor $\pi$ into a purely inseperable morphism $h:F\rightarrow F'$ and a seperable morphism $g:F'\rightarrow D$  with $F'$ a smooth projective curve. Theorem 6.3 shows that $h$ is an integer power of Frobenius and Theorem 6.5 shows that the notions of Zariski multiplicity and algebraic multiplicity coincide for the morphism $g$. Now the result follows by the fact that $h$ has algebraic multiplicity $p^{n}$ everywhere but is Zariski unramified.\\

\begin{theorem}

Let hypotheses be as in Theorem 6.5, with the modification that
$F$ and $D$ are arbitrary smooth projective varieties. Then the
notions of Zariski multiplicity and algebraic multiplicity
coincide.

\end{theorem}

We will make the necessary modifications to Theorem 6.5;\\

$(1)$. We use the following classical fact (Projective
Normalisation), see \cite{Mum}.

\begin{fact}

Let $D\subset P^{n}(L)$ be an $r$ dimensional projective variety.
Then there exist $(r+1)$ linear forms $l_{0}(X),\ldots,l_{r}(X)$
with coeffients in $L$ such that the hyperplane $H$ defined by
$l_{0}=\ldots=l_{r}=0$ is dijoint from $D$. If
$\tau:P^{n}(L)-H\rightarrow P^{r}(L)$ denotes the projection, then
the restriction of $\tau$ to $D$ is a finite surjective morphism.

\end{fact}

Now combining this with the result in Section 5 that algebraic multiplicity is multiplicative for morphisms between smooth projective varieties, we need only consider the branched finite cover $\pi:F\rightarrow P^{r}(L)$\\

$(2)$. In this case, there is no straightforward way to present
the pullbacks of the local uniformisers $x_{1},\ldots,x_{r}$ at
$\bar 0\in A^{r}$. Instead, we have the inclusion $\pi^{*}O_{\bar
0, A^{r}}\subset O_{a,F}$ induced by the map $\pi$. Passing to the
Henselisations,
 gives an inclusion ${O_{\bar O, A^{r}}^{\wedge} \subset O_{a,F}^{\wedge}}$ for the etale topology.
  In the case when $\pi$ fails to be etale in an open neighborhood of $a$, this is in fact a
  proper inclusion. Now choose uniformisers $w_{1},\ldots w_{n}$ for $O_{a,F}$. As $a$ and $\bar 0$
   are smooth points, the completions of the local rings $O_{\bar O, A^{n}}$ and $O_{a,F}$ with
    respect to the order valuations at $a$ and $\bar 0$ are isomorphic to the formal power series
     rings $L[[w_{1},\ldots, w_{n}]]$ and $L[[x_{1},\ldots,x_{n}]]$ respectively. The following
     is a classical result used in the proof of the Artin approximation theorem, relating the
      Henselisation of the ring $L\{x_{1},\ldots,x_{n}\}$ of strictly convergent power series in several
      variables with its formal completion $L[[x_{1},\ldots,x_{n}]]$.
       see \cite{Art} or \cite{Rob};\\

Henselisation$(L[x_{1},\ldots x_{n}]_{(x_{1},\ldots x_{n})})=L[[x_{1},\ldots x_{n}]]\cap L(x_{1},\ldots x_{n})^{alg}$\\

This implies that\\

$O_{\bar 0, A^{n}}^{\wedge}\cong L[[x_{1},\ldots x_{n}]]\cap L(x_{1},\ldots x_{n})^{alg}$\\

$O_{a,F}^{\wedge}\cong L[[w_{1},\ldots,w_{n}]]\cap L(w_{1},\ldots,w_{n})^{alg}$ (*)\\

We now use analytic results for the formal power series ring $L[[w_{1},\ldots, w_{n}]]$. By Weierstrass preparation, we obtain the equations\\

$x_{1}=u_{1}(w_{1}^{m_{1}}+q_{11}(w_{2},\ldots,w_{n})w_{1}^{m_{1}-1}+\ldots q_{m_{1}1}(w_{2}\ldots w_{n}))$\\

$x_{2}=u_{2}(w_{1}^{m_{2}}+q_{12}(w_{2},\ldots,w_{n})w_{1}^{m_{2}-1}+\ldots q_{m_{2}2}(w_{2}\ldots w_{n}))$\\

$\ldots$\\

$x_{n}=u_{n}(w_{1}^{m_{n}}+q_{1n}(w_{2}\ldots w_{n})w_{1}^{m_{n}-1}+\ldots q_{m_{n}n}(w_{2}\ldots w_{n}))$       (**)\\

where the $u_{i}$ are units in $L[[w_{1},\ldots,w_{n}]]$ and the $q_{ij}$ are polynomials without constant term. In order to apply Weierstrass preparation, we require that the power series expansions for the $x_{i}$ should be regular with respect to the variable $w_{1}$. Clearly this can be acheived in the following manner;\\

Let $M=(m_{kl})_{1\leq k,l\leq n}$ be an invertible matrix of elements in $L$. Then if $\bar w'=M(\bar w)$, as $M$ is invertible, $\bar w'$ is also a set of uniformisers for $L[[w_{1},\ldots w_{n}]]$. The condition of irregularity $C_{ij}$ for $x_{i}$ in terms of the variable $w_{j}'$ is a (possibly infinite) conjunction of closed relations on the $m_{kl}$. Hence there exists a Zariski open set $U\subset GL_{n}(L)$ such that $C_{ij}$ fails to hold for $1\leq i,j\leq n$, that is, after a linear change of variables,  we can assume that the $x_{i}$ each have regular expansions in terms of $w_{j}$. \\

Now $w_{1},\ldots w_{n}$ are algebraically independent in $L(F)$ which has transcendence degree $n$ over $L$. As each $x_{i}\in L(F)$, we must have that each $x_{i}\in L(w_{1},\ldots,w_{n})^{alg}$. Therefore, the $u_{i}$ in the equations $(**)$ can be taken in $L(w_{1}\ldots,w_{n})^{alg}$ and, using (*), the equations hold in $O_{a,F}^{\wedge}$.\\

(3)  Hence, we can find an etale morphism $\pi:(U,b)\rightarrow (F,a)$ such that the equations $(**)$ hold in the local ring $O_{U,b}$. Again, we may assume that $U$ is irreducible and smooth.\\

(4)-(7) This part of the argument goes through essentially unchanged, with the slight modification that the projective closure $F'$ of $U$ may fail to be smooth and the closure of $graph(\pi)$ in $F'\times F$ may fail to define a function, only a generically finite correspondence between $F'$ and $F$. However, this still allows us to work in the context of Theorem 3.3, see also the Remarks 3.4, when we consider the projection of the correspondence restricted to $U$.\\

(8) Now we have the sequence of ring inclusions given by\\

\ \ \ \ \ \ \ \ \ \ \   \ \   \ \ \ \ \ \  \ \ \ \ \ \ $L[x_{1},\ldots,x_{n}]\rightarrow\\
 L[\bar x,\bar w, \bar u]/<x_{1}-u_{1}p_{1}(\bar w),\ldots,x_{n}- u_{n}p_{n}(\bar w), s_{1}(u_{1}),\ldots s_{n}(u_{n})>\rightarrow R$\\

where $R$ is the coordinate ring of $U$ in some affine
neighborhood of $b$, $p_{i}$ are the polynomials given in $(**)$
and $s_{i}$ are the minimum polynomials of $u_{i}$ over
$L(w_{1},\ldots,w_{n})$. A simple calculation shows that the
second variety is smooth at $\bar 0$ and the second inclusion
corresponds to an etale extension of algebras. It is therefore
sufficient to check that the algebraic and Zariski multiplicities
of the left hand inclusion coincide at $\bar 0$(***). An easy
calculation gives that the algebraic multiplicity of the left hand
inclusion is $length(L[\bar w,\bar u]_{\bar 0}/<u_{1}p_{1}(\bar
w),\ldots u_{n}p_{n}(\bar w), s_{1}(u_{1}),\ldots, s_{n}(u_{n})>)$
which, by the localisation at $\bar 0$, is just $length(L[\bar
w]_{\bar 0}/<p_{1}(\bar w),\ldots p_{n}(\bar w)>)$. This is
precisely the intersection multiplicity of the hypersurfaces
$p_{1},\ldots,p_{n}$ at $\bar 0$. Again, for ease of exposition,
we compute the case of $2$ irreducible intersecting polynomials
 $p_{1}(x,y)=0$ and $p_{2}(x,y)=0$ with $p_{1}(0,0)=p_{2}(0,0)=0$. We claim the following theorem;

\begin{theorem}

 The intersection multiplicity of $p_{1},p_{2}$ at $(0,0)$ corresponds to the Zariski multiplicity
 of the cover $Spec(L[xyuv]/<p_{1}-u,p_{2}-v>)\rightarrow Spec(L[uv])$, when $Char(L)=0$.

\end{theorem}

The theorem includes the proof of $(***)$ when $n=2$. We shall
indicate how the higher
dimensional case follows later. In order to prove the theorem, we need a series of lemmas.\\

\begin{lemma}

Let $F(x,\bar y)$ be an irreducible Weierstrass polynomial in $x$
with $F(0,\bar 0)=0$ then algebraic multiplicity and Zariski
multiplicity coincide for the cover $Spec(L[x\bar
y]/<F>)\rightarrow Spec(L[\bar y])$.

\end{lemma}
\begin{proof}
We have that $F(x,\bar y)=x^{n}+q_{1}(\bar y)x^{n-1}+\ldots+q_{n}(\bar y)$ where $q_{i}(\bar 0)=0$. The algebraic multiplicity is given by $length(L[x]/F(x,\bar 0))=ord(F(x,\bar 0)=n$ in the ring $L[x]$ with the canonical valuation. We first claim that the Zariski multiplicity is the number of solutions to $x^{n}+q_{1}(\bar{\epsilon})x^{n-1}+\ldots+q_{n}(\bar{\epsilon})=0$ (*), where $\bar \epsilon$ is generic in ${\mathcal V}_{\bar 0}$. For suppose that $(a,\bar{\epsilon})$ is such a solution, then $F(a,\bar{\epsilon})=0$ and by specialisation $F(\pi(a),\bar 0)=0$. As $F$ is a Weierstrass polynomial in $x$, $\pi(a)=0$, hence $a\in{\mathcal V}_{0}$, giving the claim. As $char(L)=0$, $Disc(F(x,\bar y))=Res_{\bar y}(F,{\partial F\over \partial x})$ is a regular polynomial in $\bar y$ defined over $L$. By genericity of $\bar{\epsilon}$, we have that $Disc(F(x,\bar y))|\bar{\epsilon}\neq 0$, hence (*) has no repeated roots. This gives the lemma.\\
\end{proof}
\begin{lemma}

Let $F(x,\bar y)$ be an irreducible polynomial with $F(x,\bar
0)\neq 0$ and $F(0,\bar 0)=0$. Then the Zariski multiplicity of
the cover $Spec(L[x,\bar y]/<F>)\rightarrow Spec(L[\bar y])$
equals $ord(F(x,\bar 0))$ in $L[x]$.

\end{lemma}
\begin{proof}
By the Weierstrass Preparation Theorem, we can write $F(x,\bar y)=U(x,\bar y)G(x,\bar y)$ with $U(x,\bar y),G(x,\bar y)\in L[[x,\bar y]]$, $G(x,\bar y)$ a Weierstrass polynomial in $x$ and $deg(G)=ord(F(x,\bar 0))$. As above, we may take the new coefficients to lie inside the Henselized ring $L[x,\bar y]_{\bar 0}^{\wedge}$, hence inside some finite etale extension $L[x,\bar y]^{ext}$ of $L[x,y]$ (possibly after localising $L[x,\bar y]$). Now we have the sequence of morphisms;\\

$Sp(L[x,\bar y]^{ext}/UG)\rightarrow Spec(L[x,\bar y]/F)\rightarrow Spec(L[\bar y])$\\

The left hand morphism is etale at $\bar 0$, hence as we have seen, to compute the Zariski multiplicity of the right hand morphism, we need to compute the Zariski multiplicity of the cover\\

$Spec(L[x,\bar y]^{ext}/UG)\rightarrow Spec(L[\bar y])$\\

Choose $\bar \epsilon\in {\mathcal V}_{\bar 0}$, the fibre of the cover is given formally analytically by $L[[x,\bar y]]/<UG>\otimes_{L[\bar y],\bar y\mapsto\bar\epsilon}L$, hence by solutions to $U(x,\bar \epsilon)G(x,\epsilon)$. By definition of Zariski multiplicity, we consider only solutions $(x\bar\epsilon)$ in ${\mathcal V}_{(0,\bar 0})^{lift}$, (here $(0,\bar 0)^{lift}$ is the lift of $(0,\bar 0)$ in the etale neighborhood, for ease of notation we will just use $(0,\bar 0)$ from now on.)
As $U(x,\bar y)$ is a unit in the local ring $L[x,\bar y]^{ext}_{0,\bar 0}$, we must have $U(x,\bar\epsilon)\neq 0$ for such solutions. Hence, the solutions are given by $G(x,\bar \epsilon)=0$. Now, we use the previous lemma to give that the Zariski multiplicity is exactly $deg(G)$ as required.\\
\end{proof}
\begin{lemma}

Let $p_{1}(x,y)$, $p_{2}(x,y)$ be Weierstrass polynomials in $x$
with $p_{1}(0,0)=p_{2}(0,0)=0$. Then the Zariski multiplicity of
the cover $Spec(L[x,y,u,v]/<p_{1}-u,p_{2}-v>)\rightarrow
Spec(L[u,v])$ $(*)$ at $(\bar 0)$ equals the intersection
multiplicity at $(0,0)$, $I(p_{1},p_{2},(0,0))$.

\end{lemma}

\begin{proof}

Let $F(y,u,v)=Res(p_{1}-u,p_{2}-v)$. Then
$F(0,0,0)=Res(p_{1},p_{2})(0)=0$, as $p_{1},p_{2}$ have a common
root at $(0,0)$. By a result due to Abhyankar, see for example
\cite{Aby}, $ord_{y}(F(y,\bar 0))=\Sigma I(p_{1},p_{2},(x0))$ at
common solutions $(x,0)$ to $p_{1}$ and $p_{2}$ over $0$. As
$p_{1}$ and $p_{2}$ are Weierstrass polynomials in $x$, this is
just $I(p_{1},p_{2},(00))$. By the previous lemma, it is therefore
sufficient to prove that the Zariski multiplicity of the cover
$(*)$ at $(0,0,0,0)$ equals the Zariski multiplicity of the cover
$Spec(K[y,u,v]/<F>)\rightarrow Spec(K[u,v])$ $(**)$at $(0,0,0)$.
Suppose the Zariski multiplicity of $(**)$ equals $n$. Then there
exist $y_{1},\ldots,y_{n}\in{\mathcal V}_{0}$ distinct and $\bar
\epsilon\in{\mathcal V}_{00}$ such that $F(y_{i},\bar\epsilon)$
holds. Consider $Q(u,v)=res(F(y,u,v),{\partial F\over\partial
y}(y,u,v))$. By genericity, we have that $Q(\bar\epsilon)\neq 0$.
Hence, $F(y_{i},\bar\epsilon)$ is a non-repeated root. Using
Abhyankar's result, we can find a unique $x_{i}$ with $(x_{i}y_{i})$ a common solution to $p_{1}-\epsilon_{1}$ and $p_{2}-\epsilon_{2}$. We claim that each $(x_{i}y_{i})\in {\mathcal V}_{00}$. As $p_{1}(x_{i}y_{i})-\epsilon_{1}=0$, by specialiation $p_{1}(\pi(x_{i}),0)=0$. Now, using the fact that $p_{1}$ is a Weierstrass polynomial in $x$, gives that $\pi(x_{i})=0$ as well. This shows that the Zariski multiplicity of the cover $(*)$ is at least $n$. A virtually identical argument shows that the Zariski multiplicity of the cover $(*)$ is at most $n$ as well. This gives the result.\\

\end{proof}

\begin{lemma}

Let $p_{1}(x,y), p_{2}(x,y)$ be polynomials with
$p_{1}(0,0)=p_{2}(0,0)=0$. Then the Zariski multiplicity of the
cover $Spec(L[xyuv]/<p_{1}-u,p_{2}-v>)\rightarrow Spec(L[uv])$
equals $I(p_{1},p_{2},(00))$.

\end{lemma}

\begin{proof}
Again, using the Weierstrass Preparation Theorem, write $p_{1}(x,y)=u_{1}(x,y)f_{1}(x,y)$ and $p_{2}(x,y)=u_{2}(x,y)f_{2}(x,y)$, with $f_{1},f_{2}$ Weierstrass polynomials in $x$. As before, we may assume the new coeffiecients lie in a finite ring extension $L[x,y]^{ext}$ such that the map\\

$Spec(L[x,y]^{ext}[u,v]/<u_{1}f_{1}(x,y)-u,u_{2}f_{2}(x,y)-v>)\rightarrow Spec(L[xyuv]/<f_{1}-u,f_{2}-v>)$\\

is etale near $\bar 0$\\

Again, it is sufficient to prove that the Zariski multiplicity of the cover $Spec(L[x,y]^{ext}[u,v]/<u_{1}f_{1}(x,y)-u,u_{2}f_{2}(x,y)-v>)\rightarrow Spec(L[u,v])$ at $(0,0,0,0)$ equals $I(u_{1}f_{1},u_{2}f_{2},00)=I(f_{1},f_{2},00)$. For this we need the following ``unit removal'' lemma.\\

\begin{lemma}(Unit Removal)
\\

Let $u_{1}(x,y),u_{2}(x,y),f_{1}(x,y),f_{2}(x,y)$ be polynomials
in $L[x,y]$ with $u_{1},u_{2}$ units in the local ring
$L[x,y]_{0,0}$. Then the Zariski multiplicity of the cover
$Spec(L[x,y,u,v]/<u_{1}f_{1}(x,y)-u,u_{2}f_{2}(x,y)-v>)\rightarrow
Spec(L[u,v])$ $(*)$ is equal to the Zariski multiplicity of the
cover $Spec(L[x,y,u,v]/<f_{1}(x,y)-u,f_{2}(x,y)-v>)\rightarrow
Spec(L[u,v])$ $(**)$.

\end{lemma}

In order to prove the lemma, we first need to introduce a new version of Zariski multiplicity. Suppose that $F\subset D\times V^{n}$ is a finite cover of a smooth $2$-dimensional base $D$. \\

\begin{defn}

Given $(a,\lambda_{1},\lambda_{2})\in F$, we define;\\

 $Left.Mult_{a,\lambda_{1},\lambda_{2}}(F/D)=Card({\mathcal V}_{a}\cap F(x,\lambda_{1}',\lambda_{2}))$
 for $\lambda_{1}'\in{\mathcal V}_{\lambda_{1}}$ generic over $L$.\\
 $Right.Mult_{a,\lambda_{1},\lambda_{2}}(F/D)=Card({\mathcal
 V}_{a}\cap F(x,\lambda_{1},\lambda_{2}'))$ for
 $\lambda_{2}'\in{\mathcal V}_{\lambda_{2}}$ generic over $L$.\\
\end{defn}

By factoring the specialisations involved, it is easily shown that both left multiplicity, right multiplicity are well defined and moreover the following holds;\\

$Mult_{(a,\lambda_{1},\lambda_{2})}(F/D)=\Sigma_{a'\in({\mathcal V}_{a}\cap F(x,\lambda_{1}',\lambda_{2}))}
Right.Mult_{(a'\lambda_{1}'\lambda_{2})}(F/D)$\\

$Mult_{(a,\lambda_{1},\lambda_{2})}(F/D)=\Sigma_{a'\in({\mathcal
V}_{a}\cap
F(x,\lambda_{1},\lambda_{2}'))}Left.Mult_{(a'\lambda_{1},\lambda_{2}')}(F/D)$\\

That is we may compute the Zariski multiplicty by varying the
family in $2$ stages. Now, in the case of the lemma, after varying
one parameter, an easy algebraic calculation shows the resulting
curves intersect transversally at simple points $(x_{i}y_{i})$. In
this case we can apply the inverse function theorem to one curve
$C_{1}$ given by $u_{1}f_{1}=0$ and obtain formally analytic
presentations around each $(x_{i}y_{i})$ in the variable $t_{i}$.
As we have already seen in the previous use of analytic methods,
this does not effect the calculation of Zariski multiplicity. If
$(t_{i},h(t_{i}))$ with $h(t_{i})\in L[[t_{i}]]$ is a local
analytic presentation of $C_{1}$ at $(x_{i}y_{i})$, then, by
transversality, we have
$ord_{t_{i}}(u_{2}f_{2}(t_{i},h(t_{i})))=1$ and we have to check
that this agrees with the Zariski Right multiplicity. This
calculation has already been done in Theorem 6.5. Hence, we can
calculate the Zariski multiplicity of $(*)$ and $(**)$ as the
Zariski Left multiplicity. Now, we claim that the Zariski Left
multiplicity of the covers $(*)$ and $(**)$ is the same. This is a
straightforward calculation, suppose that the Zariski Left
Multiplicity of $(*)$ is $n$. Then there exists $\epsilon$ generic
and $(x_{1}y_{1}),\ldots,(x_{n}y_{n})\in{\mathcal V}_{00}$ such
that $u_{1}q_{1}(x_{i}y_{i})=\epsilon$ and
$u_{2}q_{2}(x_{i}y_{i})=0$. Now using the fact that the $u_{i}$
are units, we find $\epsilon'$ generic in ${\mathcal V_{0}}$ such
that $q_{1}(x_{i}y_{i})=\epsilon'$ and $q_{2}(x_{i}y_{i})=0$. This
shows
exactly that the Zariski Left Multiplicity of $(**)$ at $(0,0,0,0)$ is at least $n$.
Reversing the argument shows the Zariski Left Multiplicity is exactly $n$ as required\\

Now the proof of Lemma 6.13 follows from the proof of Lemma 6.12.
\end{proof}

Higher dimensional case; The same method as for curves, inductive
argument using Abhyankar's Lemma on resultants and Weierstrass Preparation for the ring
$L[[x_{1},\ldots,x_{n}x_{n+1}]]$.\\

\end{section}

\begin{section}{Further Directions of Study}

\begin{rmk}

Let hypotheses be as in Theorem 3.3, with the additional
assumption that $F$ is an etale cover of $D=A^{n}$. Then one can
improve the lifting condition to points in $L[[t]]$. Use the local
uniformisers to present the cover over $A^{n}$ in the form
$f_{1}(\bar x,\bar y)=0, f_{2}(\bar x,\bar y)=0,\ldots, f_{n}(\bar
x,\bar y)=0$, with $\bar x,\bar y$ tuples in $A^{n}$. Then letting
$\bar x$ be a point in $L[[t]]$ gives $n$ equations inside $A^{n}$
with coefficients in $L[[t]]$. By Hensel's lemma, we can find a
solution to these equations in $L[[t]]$, as reducing the equations
modulo $(t)$, by the fact that the morphism is etale at $\bar a$,
$\overline{f_{1}},\ldots,\overline{f_{n}}$ have a common solution
$\bar a$ in $L$ with $({\partial \overline{f_{i}}\over
\partial y_{j}})_{ij}(\bar a)\neq 0$\\

In deformation theory arguments, we work with schemes defined over the projective limit of rings
 $L[t]/(t^{n})$. This suggests developing part of the theory of Zariski structures in the analytic
  context of complete valued fields, possibly using the Pas language with sorts for the reductions
  modulo $t^{n}$. We save this point of view for another occasion. \\

\end{rmk}

\begin{rmk}

As mentioned before, one can define the etale topology and obtain
the Cech cohomology groups with finite coefficients for any
$1$-dimensional Zariski strucure. The following is a classical
result, most famously used in Deligne's proof of the Weil conjectures;\\

(Lefschetz fixed-point formula)\\

Let $X$ be a complete non-singular variety over an algebraically
closed field $K$, and let $\phi:X\rightarrow X$ be a regular map.
Then\\

$(\Gamma_{\phi}\centerdot\Delta)=\Sigma(-1)^{r}Tr(\phi|H^{r}(X,{\mathcal
Q}_{l})$  $(*)$\\

where $\Gamma_{\phi}$ is the graph of $\phi$, $\Delta$ is the
diagonal in $X\times X$ and $(\Gamma_{\phi}\centerdot\Delta)$ is
the number of fixed points of $\phi$ counted with multiplicity. \\

Now both sides of the above formula make sense in the more
generalised setting of $X$, a closed presmooth subset of $C^{n}$,
where $C$ is a $1$-dimensional Zariski structure. We use the
notion of Zariski multiplicity to replace algebraic
multiplicity. The natural question is the following;\\

For what class of Zariski structures does equality hold in
$(*)$?\\

In the algebraic context, the Lefschetz formula is a formal
consequence of a cohomology theory with good properties;\\

1. Kunneth Formula.\\
2. Finite dimensionality of the groups $H^{i}(X,F_{l^{n}})$ for
prime $l$.\\
3. Poincare duality.\\
4. Existence of a cycle map $cl^{*}_{X}:CH^{*}(X)\rightarrow H^{*}(X)$\\
5. Smooth and Proper Base Change Theorems.\\

(Here $CH^{*}(X)$ is the graded Chow ring of cycles on $X$ and
$H^{*}(X)$ is the graded cohomology ring on $X$.)\\

Clearly, an answer to the above can be reduced to further
questions concerning the class of Zariski structures for which the
 properties 1-5 hold. The interested reader should look at \cite{Milne2} or \cite{FrKi}\\

\end{rmk}

\end{section}

\end{document}